\documentclass[12pt, leqno]{article}
\usepackage{amsfonts}
\usepackage{amsmath}
\usepackage{amssymb}
\usepackage{amscd}

\title{\bf On the universal deformations for ${\rm SL}_2$-representations of knot groups}
\author{Masanori Morishita, Yu Takakura, Yuji Terashima and Jun Ueki}
\date{}

\begin{document}
\maketitle
\begin{center}
{\large \em Dedicated to Professor Kunio Murasugi}
\end{center}

\footnote[0]{
$\;\;\;\;$ Key words and phrases: Deformation of a representation, Character scheme, Knot group, Arithmetic topology \\
$\;\;\;\;$ The first author is partly supported by Grant-in-Aid for Scientific Research (B) 24340005, Japan Society for the Promotion of Science.\\
$\;\;\;\;$ The third author is partly supported by Grant-in-Aid for Scientific Research (C)  25400083, Japan Society for the Promotion of Science.\\
$\;\;\;\;$  The fourth author is partly supported by the Grant-in-Aid for JSPS Fellows (25-2241), The Ministry of Education, Culture, Sports, Science and Technology, Japan.
}

{\small{\bf Abstract.}  Based on the analogies between knot theory and number theory, we study a deformation theory for ${\rm SL}_2$-representations of knot groups, following after Mazur's deformation theory of Galois representations. Firstly, by employing the pseudo-${\rm SL}_2$-representations, we prove the existence of the universal deformation of a given ${\rm SL}_2$-representation of a finitely generated group $\Pi$ over a perfect field $k$ whose characteristic is not 2. We then show its connection with the character scheme for ${\rm SL}_2$-representations of $\Pi$ when $k$ is an algebraically closed field. We investigate examples concerning Riley  representations of 2-bridge knot groups and give  explicit forms of the universal deformations. Finally we discuss the universal deformation of the holonomy representation of a hyperbolic knot group in connection with Thurston's theory on deformations of hyperbolic structures.}
\vspace{.2cm}\\
\begin{center}
{\bf Introduction}
\end{center}

The motivation of this paper is coming from the analogies between knot theory and number theory. The study of those analogies is now  called arithmetic topology ([11]). In particular, it has been known that there are close analogies between Alexander-Fox theory and Iwasawa theory, where the Alexander polynomial and the Iwasawa polynomial ($p$-adic zeta function) are analogous objects, for instance ([8], [11, Chapters 8 -- 12]). As Mazur pointed out ([10], [11, Chapters 13, 14]), from the viewpoint of group representations, Alexander-Fox theory and Iwasawa theory are concerned about 1-dimensional representations of  knot and Galois groups, respectively, and it would be interesting to pursue the analogies further for higher dimensional representations. \\

As a first step to explore this perspective, in this paper, we study a deformation theory for representations of knot groups, following after the deformation theory for  Galois representations ([9]). In fact, we develop a general theory on deformations for ${\rm SL}_2$-representations of a finitely generated group. We deal with only ${\rm SL}_2$-representations, since our main interest is applications to knot theory and 3-dimensional topology (hyperbolic geometry) where the character varieties of ${\rm SL}_2$-representations of fundamental groups have often been studied. See [3]  for example. We note that while Galois deformation theory is concerned with $p$-adic deformation of a continuous representation of a profinite group over a finite residue field, our work deals with  infinitesimal deformation of a representation of any finitely generated group over any   perfect residue field whose characteristic is not 2 (for example, the field of complex numbers), and that our universal deformation space may be regarded as an infinitesimal (1-parameter) deformation over a complete local algebra of the character variety over a filed (see Theorem 3.2.1). Thus our work is applicable to geometry and topology. \\

The contents of this paper are as follows. In Section 1, following Wiles ([23] for ${\rm GL}_2$ case) and Taylor ([20] for ${\rm GL}_n$ case), we introduce the notion of a pseudo-${\rm SL}_2$-representation of $\Pi$ over a commutative ring and prove the existence of the universal deformation of a given pseudo-${\rm SL}_2$-representation over a perfect field.  In Section 2,  for a given representation over a perfect field $k$ whose characteristic is not 2

$$\overline{\rho} : \Pi \longrightarrow {\rm SL}_2(k),$$
 we prove, using the result in Section 1,  that there exists the universal deformation of $\overline{\rho}$
$$ {\boldsymbol \rho} : \Pi \longrightarrow {\rm SL}_2({\boldsymbol R}_{\overline{\rho}}),  $$ 
which parametrizes all lifts of $\overline{\rho}$ to ${\rm SL}_2$-representations over complete local ${\cal O}$-algebras where ${\cal O}$ a complete discrete valuation ring whose residue field is $k$.
A merit to make use of pseudo-representations is to enable us to relate the universal deformation ring with the character scheme/variety of ${\rm SL}_2$-representations where the latter has been extensively studied in the context of topology (e.g., [3], [6], [16] etc).  In fact, in Section 3,  when $k$ is an algebraically closed field, we show the relation between the universal deformation ring ${\boldsymbol R}_{\overline{\rho}}$ and the ${\rm SL}_2$-character scheme  of $\Pi$.
In Section 4, we investigate examples concerning Riley representations of 2-bridge knot groups ([17]) and give explicit forms of universal deformations. In Section 5, we apply our deformation theory to the case where $\Pi$ is the fundamental group of the complement of a hyperbolic knot in the 3-sphere and $\overline{\rho}$ is the associated holonomy representation, and describe the universal deformation ring by Thurston's deformation theory of hyperbolic structures ([21]).  We observe that our result is similar to the case of $p$-adic ordinary Galois representations where the universal deformation is described by Hida's deformation of $p$-adic ordinary modular forms ([4], [5]).\\
\\
ACKNOWLEDGMENTS. We would like to thank Gregory Brumfiel, Shinya Harada, Takahiro Kitayama, Tomoki Mihara, Sachiko Ohtani, Adam Sikora, and Seidai Yasuda for useful communications.  \\
\\
{\bf Notation.}  For a local ring $R$, we denote by $\frak{m}_R$ the maximal ideal of $R$. For an integral domain $k$, we denote by ${\rm char}(k)$ the characteristic of $k$.\\

\begin{center}
{\bf 1. Pseudo-representations and their deformations}
\end{center}

In Sections 1 and 2, we develop a deformation theory of representations for any finitely generated group. We consider only  ${\rm SL}_2$-representations, since our main concern is applications to knot theory and 3-dimensional topology where ${\rm SL}_2$-representations of  fundamental groups have been often studied. See [3] for example. Moreover, while Galois deformation theory is concerned with $p$-adic deformation of a continuous representation of a profinite group over a finite residue field, we study infinitesimal deformation of a representation of any finitely generated group over any perfect residue field whose characteristic is not 2 (for example, the field of complex numbers). Thus our work is applicable to geometry and topology.\\

In Subsection 1.1, we introduce the notion of a pseudo-${\rm SL}_2$-representation of a finitely generated group. This notion was originally introduced by Wiles ([23] for ${\rm GL}_2$ case) and by Taylor ([20] for ${\rm GL}_n$ case)  with the intention of applications to $p$-adic Galois representations. In Subsection 1.2,  we show the existence of the universal deformation of a given pseudo-${\rm SL}_2$-representation over any perfect field. \\

1.1. {\em Pseudo-${\rm SL}_2$-representations.} Let $\Pi$ be a finitely generated group. Let $A$ be  a commutative ring with identity. A map $T : \Pi \rightarrow A$ is called a {\em pseudo-${\rm SL}_2$-representation} over $A$ if the following four conditions are satisfied:
\vspace{.18cm}\\
(P1) $T(1) = 2$,\\
(P2) $T(g_1g_2) = T(g_2g_1)$ for any $g_1, g_2 \in \Pi$,\\
(P3) $T(g_1)T(g_2)T(g_3) + T(g_1g_2g_3) + T(g_1g_3g_2) -T(g_1g_2)T(g_3)  - T(g_2g_3)T(g_1) - T(g_1g_3)T(g_2)  = 0$ for any $g_1, g_2, g_3\in \Pi$,\\
(P4) $T(g)^2 - T(g^2) = 2$ for any $g \in \Pi$.
\vspace{.18cm}\\
Note that the conditions (P1) through (P3) are nothing but Taylor's conditions for a  pseudo-representation of degree 2 ([20]) and that (P4) is the condition for determinant 1. By the invariant theory of matrices ([15, Theorem 4.3]), the trace ${\rm tr}(\rho)$ of a representation $\rho : \Pi \rightarrow {\rm SL}_2(A)$ satisfies the conditions (P1) through (P4). Conversely, a pseudo-${\rm SL}_2$-representation is shown to be obtained as the trace of a representation under certain conditions (See Theorem 2.2.1 below). \\

1.2. {\em Deformations of pseudo-${\rm SL}_2$-representations.} We fix a perfect field $k$ and a complete discrete valuation ring ${\cal O}$ with the  residue field ${\cal O}/\frak{m}_{\cal O} = k$. We may take ${\cal O}$ to be the Witt ring of $k$ if  ${\rm char}({\cal O}) \neq {\rm char}(k)$, and ${\cal O} = k[[\hbar]]$, the formal power series ring of a variable $\hbar$ over $k$, if ${\rm char}({\cal O}) = {\rm char}(k)$.  There is a unique subgroup $V$ of ${\cal O}^{\times}$ such that $k^{\times} \simeq V$ and ${\cal O}^{\times} = V \times (1+\frak{m}_{\cal O})$. The composition map $\varphi : k^{\times} \simeq V \hookrightarrow {\cal O}^{\times}$ is called the {\em Teichm\"{u}ller lift} which satisfies $\varphi(\alpha) \, \mbox{mod}\, \frak{m}_{\cal O} = \alpha$ for $\alpha \in k$.  It is extended to $\varphi : k \hookrightarrow {\cal O}$ by $\varphi(0) := 0$. Let ${\cal C}$ be the category of complete local ${\cal O}$-algebras with residue field $k$. A morphism in ${\cal C}$ is an ${\cal O}$-algebra homomorphism inducing the identity on residue fields. 

Let $\overline{T} : \Pi \rightarrow k$ be a pseudo-${\rm SL}_2$-representation over $k$. A couple $(R, T)$ is called an {\em ${\rm SL}_2$-deformation} of $\overline{T}$ if $R \in {\cal C}$ and $T : \Pi \rightarrow R$ is a pseudo-${\rm SL}_2$-representation over $R$ such that $T$ mod $\frak{m}_R = \overline{T}$. In the following, we say simply a {\em deformation} of $\overline{T}$ for an ${\rm SL}_2$-deformation. A deformation $({\boldsymbol R}_{\overline{T}}, {\boldsymbol T})$ of $\overline{T}$ is called a {\em universal deformation} of $\overline{T}$ if the following universal property is satisfied: ``For any deformation $(R, T)$ of $\overline{T}$ there exists a unique morphism $\psi : {\boldsymbol R}_{\overline{T}} \rightarrow R$ in ${\cal C}$ such that $\psi \circ {\boldsymbol T} = T$."  So the correspondence $\psi \mapsto \psi \circ {\boldsymbol T}$ gives the bijection 
$${\rm Hom}_{\cal C}({\boldsymbol R}_{\overline{T}}, R) \simeq \{(R,T) \, | \, \mbox{deformation of} \; \overline{T} \}. $$
Note that a universal deformation of $\overline{T}$ is unique (if it exists) up to ${\cal O}$-isomorphism in the obvious sense. The ${\cal O}$-algebra ${\boldsymbol R}_{\overline{T}}$ is called the {\em universal deformation ring} of $\overline{T}$.\\
\\
{\bf Theorem 1.2.1.}  {\em For a pseudo-${\rm SL}_2$-representation $\overline{T} : \Pi \rightarrow k$, there exists a universal deformation $({\boldsymbol R}_{\overline{T}}, {\boldsymbol T})$ of $\overline{T}$.}\\
\\
{\em Proof.}  Let ${\cal R} := {\cal O}[[X_g; g \in \Pi]]$ be the ring of formal power series over ${\cal O}$ with variables $X_g$ indexed by elements of $\Pi$. By definition, the ring ${\cal R}$ consists of formal power series of variables $X_{g_i}$'s where indices $g_i$'s belong to a finite subset of $G$. Let $\varphi : k \hookrightarrow {\cal O}$ be the Teichm\"{u}ller lift. We set $T_g := X_g + \varphi(\overline{T}(g))$ for $g \in G$. Consider the ideal ${\cal I}$ of ${\cal R}$ generated by the elements of following type:
\vspace{.18cm}\\
(1) $T_1 - 2 = X_1 +\varphi(\overline{T}(1)) - 2$,\\
(2) $T_{g_1g_2} - T_{g_2g_1} = X_{g_1g_2} - X_{g_2g_1}$,\\
(3) $T_{g_1}T_{g_2}T_{g_3} + T_{g_1g_2g_3} + T_{g_1g_3g_2} -T_{g_1g_2}T_{g_3} - T_{g_2g_3}T_{g_1} - T_{g_1g_3}T_{g_2}$,\\
(4) $T_g^2 - T_{g^2} -2,$
\vspace{.18cm}\\
where $g, g_1, g_2, g_3 \in \Pi$. We then set ${\boldsymbol R}_{\overline{T}} := {\cal R}/{\cal I}$ and define a map ${\boldsymbol T} : \Pi \rightarrow {\boldsymbol R}_{\overline{T}}$ by ${\boldsymbol T}(g) := T_g$ mod ${\cal I}$. Then we note that ${\boldsymbol R}_{\overline{T}} \in {\cal C}$, and by the conditions (P1) through (P4), ${\boldsymbol T} : \Pi \rightarrow {\boldsymbol R}_{\overline{T}}$ is a pseudo-${\rm SL}_2$-representation and ${\boldsymbol T}$ mod $\frak{m}_{{\boldsymbol R}_{\overline{T}}} = \overline{T}$. Hence $({\boldsymbol R}_{\overline{T}}, {\boldsymbol T})$ is a deformation  of $\overline{T}$.

Next let $(R,T)$ be any deformation of $\overline{T}$. Define a morphism $\psi : {\cal R} \rightarrow R$ in ${\cal C}$ by $ \psi(f(X_g)) := f(T_g - \varphi(\overline{T}(g)))$ for $f(X_g) \in {\cal R}$. Note that $T_g - \varphi(\overline{T}(g)) \in \frak{m}_R$ and hence $f(T_g - \varphi(\overline{T}(g)))$ is well-defined since $R$ is complete with respect to the $\frak{m}_R$-adic topology. By (P1) through (P4), $\psi({\cal I}) = 0$ and hence we have the induced ${\cal O}$-algebra homomorphism in ${\cal C}$, denoted by the same $\psi$, $\psi : {\boldsymbol R}_{\overline{T}} \rightarrow R.$
 Then we easily see that $\psi \circ {\boldsymbol T} = T$. The uniqueness of $\psi$ follows from the fact that ${\boldsymbol R}_{\overline{T}}$ is generated by $X_g$ ($g \in \Pi$) as an ${\cal O}$-algebra. $\;\; \Box$\\
\\

\begin{center}
{\bf 2. The universal deformation for representations}
\end{center}

In this section, we are concerned with deformations of ${\rm SL}_2$-representations of a finitely generated group $\Pi$. 

In Subsection 2.1, we recall two theorems due to Carayol [2] and Nyssen [13].  In Subsection 2.2, by using them, we prove that there is a bijective correspondence given by the trace between ${\rm SL}_2$-representations (up to strict equivalence) and   pseudo-${\rm SL}_2$-representations, and then derive the existence of the universal deformation of an ${\rm SL}_2$-representation over a field. \\

2.1. {\em Carayol's and Nyssen's theorems.}  Two representations $\rho, \rho' : \Pi \rightarrow {\rm GL}_n(A)$ over a commutative ring $A$ with identity are said to be {\em equivalent}, written as $\rho \sim \rho'$, if there is $\gamma \in {\rm GL}_n(A)$ such that $\rho'(g) = \gamma^{-1} \rho(g) \gamma$ for any $g \in \Pi$. When $A$ is a local ring, $\rho, \rho'$ are said to be {\em strictly equivalent}, written as $\rho \approx \rho'$, if there is $\gamma \in I_n + {\rm M}_n(\frak{m}_A)$ such that $\rho'(g) = \gamma^{-1} \rho(g) \gamma$ for any $g \in \Pi$. We say that a representation $\rho : \Pi \rightarrow {\rm GL}_n(k)$ over a field $k$ is {\em absolutely irreducible} if for an algebraic closure $\overline{k}$ of $k$ the composite of $\rho$ with the inclusion ${\rm GL}_n(k) \hookrightarrow {\rm GL}_n(\overline{k})$ is an irreducible representation. This condition is independent of the choice of an algebraic closure $\overline{k}$. We recall the following theorem due to Carayol and Serre.\\
\\
{\bf Theorem 2.1.1} ([2, Theorem 1]). {\em Let $\rho, \rho' : \Pi \rightarrow {\rm GL}_n(A)$ be representations over a local ring $A$ with the residue field $k = A/\frak{m}_A$. If the residual representation $\rho$ mod $\frak{m}_A : \Pi \rightarrow {\rm GL}_n(k)$ is absolutely irreducible and ${\rm tr}(\rho) = {\rm tr}(\rho')$, then we have $\rho \sim \rho'$.}\\
\\
Next we recall the degree 2 case of a theorem by Nyssen. \\
\\
{\bf Theorem 2.1.2} ([13, Theorem 1]).  {\em Let $A$ be a Henselian separated local ring with the residue field $k := A/\frak{m}_A$ and let $T : \Pi \rightarrow A$ be a Taylor's pseudo-representation of degree $2$ over $A$. Assume that there is an absolutely irreducible representation $\overline{\rho} : \Pi \rightarrow {\rm GL}_2(k)$ such that ${\rm tr}(\overline{\rho}) = T$ mod $\frak{m}_A$. Then there exists a unique representation $\rho : \Pi \rightarrow {\rm GL}_2(A)$ such that ${\rm tr}(\rho) = T$.}\\

2.2. {\em Deformations of an ${\rm SL}_2$-representations.}  As in Subsection 1.2, let us fix a perfect field $k$ and a complete discrete valuation ring ${\cal O}$ with the residue field ${\cal O}/\frak{m}_{\cal O} = k$. We assume ${\rm char}(k) \neq 2$. Let ${\cal C}$ be the category of complete local ${\cal O}$-algebras with residue field $k$ where a morphism  is an ${\cal O}$-algebra homomorphism inducing the identity on residue fields. We note that $2$ is invertible in ${\cal O}$ and hence in any $R \in {\cal C}$. Let $\overline{\rho} : \Pi \rightarrow {\rm SL}_2(k)$ be a given representation. We call a couple $(R, \rho)$ an ${\rm SL}_2$-{\em deformation} of $\overline{\rho}$ if $R \in {\cal C}$ and $\rho : \Pi \rightarrow {\rm SL}_2(R)$ is a representation such that $\rho$ mod $\frak{m}_R = \overline{\rho}$. In the following, as in the case of pseudo-${\rm SL}_2$-representations, we say simply a {\em deformation} of $\overline{\rho}$ for an ${\rm SL}_2$-deformation. A deformation $({\boldsymbol R}_{\overline{\rho}}, \mbox{\boldmath $\rho$})$ of $\overline{\rho}$ is called a {\em universal deformation} of $\overline{\rho}$ if the following universal property is satisfied: ``For any deformation $(R, \rho)$ of $\overline{\rho}$ there exists a unique morphism $\psi : {\boldsymbol R}_{\overline{\rho}} \rightarrow R$ in ${\cal C}$ such that $\psi \circ \mbox{\boldmath $\rho$} \approx \rho$." So the correspondence $\psi \mapsto \psi \circ \mbox{\boldmath $\rho$}$ gives the bijection 
$$ {\rm Hom}_{\cal C}({\boldsymbol R}_{\overline{\rho}},R) \simeq \{ (R, \rho)\, | \, \mbox{deformation of}\; \overline{\rho}\}/\approx. $$
Note that a universal deformation of $\overline{\rho}$ is unique (if it exists) up to ${\cal O}$-isomorphism in the obvious sense. The ${\cal O}$-algebra ${\boldsymbol R}_{\overline{\rho}}$ is called the {\em universal deformation ring} of $\overline{\rho}$. 

A deformation  $(R, \rho)$ of $\overline{\rho}$ gives rise to a deformation $(R, {\rm tr}(\rho))$ of the pseudo-${\rm SL}_2$-representation ${\rm tr}(\overline{\rho}) : \Pi \rightarrow k$. The following theorem asserts that this correspondence is actually  bijective under the assumption that $\overline{\rho}$ is absolutely irreducible.\\
\\
{\bf Theorem 2.2.1.} {\em Let $\overline{\rho} : \Pi \rightarrow {\rm SL}_2(k)$ be an absolutely irreducible representation and let $R \in {\cal C}$. Then the correspondence $\rho \mapsto {\rm tr}(\rho)$ gives the following bijection}:
$$ \begin{array}{l}
\{ \rho : \Pi \rightarrow {\rm SL}_2(R) \, | \, \mbox{deformation of}\; \overline{\rho} \; \mbox{over}\; R\}/\approx  \\
\;\;\;\;\;\;\;\;\;\;\;\;\;\;\;\;\;\;  \longrightarrow \{ T : \Pi \rightarrow R \,  | \, \mbox{deformation of} \; {\rm tr}(\overline{\rho}) \; \mbox{over} \; R \}.
\end{array}
$$
\\
{\em Proof.} Firstly let us show the surjectivity. Let $T : \Pi \rightarrow R$ be a pseudo-${\rm SL}_2$-representation such that $T$ mod $\frak{m}_R = {\rm tr}(\overline{\rho})$. By Theorem 2.1.2, there exists a unique representation $\rho_1 : \Pi \rightarrow {\rm GL}_2(R)$ such that ${\rm tr}(\rho_1) = T$. Note that $\rho_1$ is actually an ${\rm SL}_2(R)$-representation, because we have $2 \det(\rho_1(g)) = {\rm tr}(\rho_1(g))^2 - {\rm tr}(\rho_1(g^2)) = T(g)^2 - T(g^2) = 2$ for any $g \in \Pi$ and $2 \in R^{\times}$. Since ${\rm tr}(\rho_1 \; {\rm mod}\; \frak{m}_R) = T \; {\rm mod}\; \frak{m}_R = {\rm tr}(\overline{\rho})$ and $\overline{\rho}$ is absolutely irreducible, Theorem 2.1.1 implies that $\overline{\rho} \sim \rho_1$ mod $\frak{m}_R$. So there is $\overline{\gamma} \in {\rm GL}_2(k)$ such that 
$\overline{\rho}(g) = \overline{\gamma}^{-1} (\rho_1\;{\rm mod} \; \frak{m}_R)(g) \overline{\gamma}$. Choose a lift $\gamma \in {\rm GL}_2(R)$ of $\overline{\gamma}$ and define a representation $\rho : \Pi \rightarrow {\rm SL}_2(R)$ by $\rho(g) := \gamma^{-1}\rho_1(g) \gamma$ for $g \in \Pi$. Then $(R, \rho)$ is a deformation of $\overline{\rho}$ and ${\rm tr}(\rho) = {\rm tr}(\rho_1) = T$.

Next let us show the injectivity. Let $\rho, \rho' : \Pi \rightarrow {\rm SL}_2(R)$ be deformations of $\overline{\rho}$ such that ${\rm tr}(\rho) = {\rm tr}(\rho')$. Since $\overline{\rho}$ is absolutely irreducible, Theorem 2.1.1 implies $\rho \sim \rho'$. So there is $\gamma \in {\rm GL}_2(R)$ such that $\rho'(g) = \gamma^{-1}\rho(g)\gamma$ for $g \in \Pi$. Taking mod $\frak{m}_R$, we have $\overline{\rho}(g) = \overline{\gamma}^{-1} \overline{\rho}(g) \overline{\gamma}$ for $g \in \Pi$ where we put $\overline{\gamma} := \gamma$ mod $\frak{m}_R$. Since $\overline{\rho}$ is irreducible, Schur's lemma implies that $\overline{\gamma}$ is a scaler matrix over $k$, say $\overline{\gamma} = \overline{a} I_2$, $\overline{a} \in k^{\times}$. Take a lift $a \in R^{\times}$ of $\overline{a}$ and set $\gamma' := aI_2$. Then $\gamma \gamma'^{-1} \equiv I_2$ mod $\frak{m}_R$ and $\rho'(g) = (\gamma  \gamma'^{-1})^{-1} \rho(g) (\gamma \gamma'^{-1})$ for $g \in \Pi$. Hence $\rho' \approx \rho$. $\;\; \Box$\\
\\
{\bf Theorem 2.2.2.} {\em Let $\overline{\rho} : \Pi \rightarrow {\rm SL}_2(k)$ be an absolutely irreducible representation. Then there exists the universal deformation $({\boldsymbol R}_{\overline{\rho}}, {\boldsymbol \rho})$ of $\overline{\rho}$, where ${\boldsymbol R}_{\overline{\rho}}$ is given as ${\boldsymbol R}_{\overline{T}}$ for $\overline{T}:= {\rm tr}(\overline{\rho})$ in Theorem 1.2.1.}\\
\\
{\em Proof.}  By Theorem 1.2.1, there exists the universal deformation $({\boldsymbol R}_{\overline{T}}, {\boldsymbol T})$ of a pseudo-${\rm SL}_2$-representation $\overline{T} = {\rm tr}(\overline{\rho})$.  By Theorem 2.2.1, we have a deformation $\mbox{\boldmath $\rho$} : \Pi \rightarrow {\rm SL}_2({\boldsymbol R}_{\overline{T}})$ of $\overline{\rho}$ such that ${\rm tr}(\mbox{\boldmath $\rho$}) = {\boldsymbol T}$. We claim that $({\boldsymbol R}_{\overline{T}}, \mbox{\boldmath $\rho$})$ is the universal deformation of $\overline{\rho}$ and hence ${\boldsymbol R}_{\overline{\rho}} = {\boldsymbol R}_{\overline{T}}$.  Let $(R, \rho)$ be any deformation of $\overline{\rho}$. By the universality of $({\boldsymbol R}_{\overline{T}}, {\boldsymbol T})$, there exists a unique morphism $\psi : {\boldsymbol R}_{\overline{T}} \rightarrow R$ in ${\cal C}$ such that $\psi \circ {\boldsymbol T} = {\rm tr}(\rho)$. Since ${\rm tr}(\psi \circ \mbox{\boldmath $\rho$}) = \psi \circ {\rm tr}(\mbox{\boldmath $\rho$}) = \psi \circ {\boldsymbol T} = {\rm tr}(\rho)$, Theorem 2.2.1 implies $\psi \circ \mbox{\boldmath $\rho$} \approx \rho$. $\;\; \Box$\\

Finally we recall a basic fact on a presentation of a complete local ${\cal O}$-algebra, which will be used later. For $R \in {\cal C}$, we define the {\em relative cotangent space} $\frak{t}_{R/{\cal O}}^*$ of $R$ by the $k$-vector space  $\frak{m}_{R}/(\frak{m}_R^2 + \frak{m}_{\cal O}R)$ and the {\em relative tangent space} $\frak{t}_{R/{\cal O}}$ of $R$ by the dual $k$-vector space of $\frak{t}_{R/{\cal O}}^*$. We note that they are same as the cotangent and tangent spaces of $R/\frak{m}_{\cal O}R = R \otimes_{\cal O} k$, respectively. The following Lemma 2.2.3 may be a well-known fact which is proved using Nakayama's lemma (cf. [22, Lemma 5.1]).\\
\\
{\bf Lemma 2.2.3.} {\em Let $ d$ be the dimension of $\frak{t}_{R/{\cal O}}$ over $k$ and assume $d < \infty$.   For a given system of parameters $x_1,\dots, x_d$ of the local $k$-algebra $R \otimes_{\cal O} k$,  there is a surjective ${\cal O}$-algebra homomorphism  $$\lambda : {\cal O}[[X_1,\dots, X_d]] \longrightarrow R$$
in ${\cal C}$ such that the image of $\lambda(X_i)$ in $R\otimes_{\cal O} k$ is $x_i$ $(1\leq i \leq d)$.}\\

\begin{center}
{\bf 3. Character schemes} 
\end{center}

In this section, we show the relation between the universal deformation ring in Sections 1, 2 and the character scheme of ${\rm SL}_2$-representations. 

In Subsection 3.1, we recall the constructions and some facts concerning the ${\rm SL}_2$-character scheme and the skein algebra over an algebraically closed field and then describe their relation.
For the details, we consult [3], [7, Chapter 1], [12] and [19]. In Subsection 3.2, we give the relation between the universal deformation ring  and the character scheme via the skein algebra. Our universal deformation ring may be regarded as an infinitesimal (1-parameter) deformation of the universal character algebra.

As in Sections 1 and 2, let $\Pi$ be a finitely generated group. 
\\

3.1. {\em Character schemes and skein algebras.} Let $k$ be an algebraically closed field and consider the functor $\frak{F}$ from the category of commutative $k$-algebras to the category of sets defined by 
$$ A \; \mapsto \; \frak{F}(A) := \mbox{the set of all representations} \; \Pi \rightarrow {\rm SL}_2(A).$$
The functor $\frak{F}$ is representaed by a commutative $k$-algebra $\frak{A}(\Pi)$, called the {\em universal ${\rm SL}_2$-representation algebra} of $\Pi$ over $k$, and we have  the {\em universal ${\rm SL}_2$-representation} $\rho^{\rm univ} : \Pi \rightarrow {\rm SL}_2(\frak{A}(\Pi))$ which satisfies the following  property:  ``For any commutative $k$-algebra $A$ and a representation $\rho : \Pi \rightarrow {\rm SL}_2(A)$, there is a unique $k$-algebra homomorphism $\psi : \frak{A}(\Pi) \rightarrow A$ such that $\rho = \psi \circ \rho^{\rm univ}$." In fact, when $\Pi$ is given by generators $g_1,\dots,g_s$ subject to the relations $r_q = 1 \;(q \in Q)$, the universal ${\rm SL}_2$-representation algebra $  \frak{A}(\Pi)$ is given by the quotient of the polynomial ring $k[X_{ij}^{(m)} \; (1\leq m \leq s, 1\leq i,j\leq 2)]$ by the ideal $J$ generated by $\det(X^{(m)}) -1 \, (1\leq m \leq s)$ and $(r_q)_{ij} \, (q \in Q, 1\leq i,j\leq 2),$
 where $X^{(m)} := (X_{ij}^{(m)})$ and $(r_q)_{ij}$ denotes the $(i,j)$-entry of the matrix $r_q(X^{(1)},\dots,X^{(s)})$, and the universal representation $\rho^{\rm univ}$ is given by $\rho^{\rm univ}(g_m) = X^{(m)}$ mod $J$ for $1\leq m \leq s$. We denote the affine scheme ${\rm Spec}(\frak{A}(\Pi) )$ by $\frak{R}(\Pi)$ and call it the {\it  ${\rm SL}_2$-representation scheme} of $\Pi$ over $k$. We identify  a prime ideal $\frak{p} \in \frak{R}(\Pi)$  with the corresponding representation $\rho_\frak{p} := \psi_\frak{p} \circ \rho^{\rm univ} : \Pi \rightarrow {\rm SL}_2(\frak{A}(\Pi) /\frak{p})$, where $\psi_\frak{p} : \frak{A}(\Pi)  \rightarrow \frak{A}(\Pi) /\frak{p}$ is the natural homomorphism. We set ${\cal R}(\Pi) := \frak{R}(\Pi)(k) = {\rm Spm}(\frak{A}(\Pi) )$ and call it  the {\it  ${\rm SL}_2(k)$-representation variety} of $\Pi$. It is an affine algebraic set over $k$ which parametrizes all representations $\Pi \rightarrow {\rm SL}_2(k)$, obtained as $\rho_{\frak{m}} = \psi_{\frak{m}} \circ \rho^{\rm univ}$ for $\frak{m} \in {\cal R}(\Pi)$, where $\psi_{\frak{m}} : \frak{A}(\Pi) \rightarrow \frak{A}(\Pi)/\frak{m} = k$ is the natural homomorphism. We identify a maximal ideal $\frak{m} \in {\cal R}(\Pi)$ and the corresponding representation $\rho_{\frak{m}}$. We denote by $k[{\cal R}(\Pi)]$ the coordinate ring of ${\cal R}(\Pi)$. We note that $k[{\cal R}(\Pi)]$ is the quotient of $\frak{A}(\Pi)$ by the nilradical, $ k[{\cal R}(\Pi)]= \frak{A}(\Pi)/\sqrt{0}$.

The adjoint action of the group scheme ${\rm GL}_2$  on $\frak{B}(\Pi)$ is defined by  sending the $(i,j)$-entry of $X^{(m)}$ to the $(i,j)$-entry of $P^{-1}X^{(m)}P$ for $P \in {\rm GL}_2$. Let $\frak{B}(\Pi)$ be the invariant subalgebra of $\frak{A}(\Pi)$ under this action of ${\rm GL}_2$, $\frak{B}(\Pi) := \frak{A}(\Pi)^{{\rm GL}_2}$, which we call  the {\em universal ${\rm SL}_2$-character algebra} of $\Pi$ over $k$. We denote by $\frak{X}(\Pi)$ the affine scheme ${\rm Spec}(\frak{B}(\Pi))$, namely, the algebro-geometric quotient of $\frak{R}(\Pi)$ by the adjoint action of ${\rm GL}_2$, and call it the {\it ${\rm SL}_2$-character scheme} of $\Pi$ over $k$. We have a morphism $\frak{R}(\Pi) \rightarrow \frak{X}(\Pi)$ induced by the inclusion $\frak{B}(\Pi) \hookrightarrow \frak{A}(\Pi)$. We write $[\frak{p}] (= [\rho_{\frak{p}}])$ for the image of $\frak{p} (= \rho_{\frak{p}})$. 
We set ${\cal X}(\Pi) := \frak{X}(\Pi)(k) = {\rm Spm}(\frak{B}(\Pi))$ and call it the {\it ${\rm SL}_2(k)$-character variety} of $\Pi$. It is an algebraic set which parametrizes all characters ${\rm tr}(\rho)$ of representations $\rho : \Pi \rightarrow {\rm SL}_2(k)$. Under the natural morphism ${\cal R}(\Pi) \rightarrow {\cal X}(\Pi)$, we write $[\rho] \in {\cal X}(\Pi)$ for the image of $\rho \in {\cal R}(\Pi)$. We note that $[\rho] = [\rho']$ if and only if ${\rm tr}(\rho) = {\rm tr}(\rho')$. We denote by $k[{\cal X}(\Pi)]$ the coordinate ring of ${\cal X}(\Pi)$. We note that $k[{\cal X}(\Pi)]$ is the invariant subring of $k[{\cal R}(\Pi)]$ under the conjugate action of ${\rm GL}_2(k)$, $k[{\cal X}(\Pi)] = k[{\cal R}(\Pi)]^{{\rm GL}_2(k)}$ and $k[{\cal X}(\Pi)] = \frak{B}(\Pi)/\sqrt{0}$.  For $g \in \Pi$, define $\tau_g : {\cal R}(\Pi) \rightarrow k$ by $\tau_g(\rho) := {\rm tr}(\rho(g))$. It is known ([LM, Corollary 1.34]) that $k[{\cal X}(\Pi)]$ is generated over $k$ by finitely many $\tau_g$'s. 
 
 According to [19, 3.1] and [16, Definition 2.5], we define  the $k$-algebra $\frak{C}(\Pi)$ by
 $$ \frak{C}(\Pi) := k[t_g \; (g \in \Pi)]/I,$$
where $t_g$ is a variable for each $g \in \Pi$ and $I$ is the ideal of the polynomial ring $k[t_g  \; (g \in \Pi)]$ generated by the polynomials of the form
$$t_1 -2,\;\;  t_{g_1}t_{g_2} - t_{g_1g_2} - t_{g_1^{-1}g_2}\;  (g_1, g_2 \in \Pi).$$
We call $\frak{C}(\Pi)$ the {\em skein algebra} of $\Pi$ over $k$. We note that $\frak{C}(\Pi)$ is Noetherian since $\Pi$ is finitely generated. We denote the affine scheme ${\rm Spec}(\frak{C}(\Pi))$ by $\frak{X}^{\rm skein}(\Pi)$ and call it the {\it skein scheme} of $\Pi$ over $k$. Since ${\rm tr}(\rho^{\rm univ}(g)) \in \frak{B}(\Pi)$ for $g \in \Pi$ and we have the formula, which is derived by the Cayley-Hamilton theorem,
$$ {\rm tr}(\rho^{\rm univ}(g_1)) {\rm tr}(\rho^{\rm univ}(g_2)) - {\rm tr}(\rho^{\rm univ}(g_1g_2)) - {\rm tr}(\rho^{\rm univ}(g_1^{-1}g_2)) = 0\;\, (g_1, g_2 \in \Pi),$$
we obtain a $k$-algebra homomorphism 
$$ \iota_{\Pi} : \frak{C}(\Pi) \longrightarrow \frak{B}(\Pi)$$
defined by $\iota(t_g) := {\rm tr}(\rho^{\rm univ}(g))$ for $g \in \Pi$, and hence a morphism of schemes
$$ \iota_{\Pi}^{a} : \frak{X}(\Pi) \longrightarrow \frak{X}^{\rm skein}(\Pi).$$
We define the {\it discriminant ideal} $\Delta(\Pi)$ of $\frak{C}(\Pi)$ by the ideal generated by the images of the elements in $k[t_g \; (g \in \pi)]$ of the form
$$ \Delta(g_1,g_2) := t_{g_1g_2g_1^{-1}g_2^{-1}} -2 = t_{g_1}^2 + t_{g_2}^2 + t_{g_1g_2}^2 - t_{g_1}t_{g_2}t_{g_1g_2} - 4 \;\; (g_1,g_2 \in \pi),$$
and the {\it discrimiant subscheme} by $V(\Delta(\Pi)) = {\rm Spec}(\frak{C}(\Pi)/\Delta(\Pi))$. Since $\frak{C}(\Pi)$ is Noetherian, $\Delta$ is generated by finitely many $\Delta(g_1^{(i)},g_2^{(i)})$, $i = 1,\dots,n$. 
We set $\Delta := \Delta(g_1^{(1)},g_2^{(2)})\cdots \Delta(g_1^{(n)},g_2^{(n)}) \in \frak{C}(\Pi)$ and define the open subschemes $\frak{X}^{\rm skein}(\Pi)_{\rm irr}$ and $\frak{X}(\Pi)_{\rm irr}$ of $\frak{X}^{\rm skein}(\Pi)$ and $\frak{X}(\Pi)$, respectively,  by
$$\begin{array}{c} \frak{X}^{\rm skein}(\Pi)_{\rm irr} := \frak{X}^{\rm skein}(\Pi) \setminus V(\Delta(\Pi)) = \frak{X}^{\rm skein}(\Pi)_{\Delta}, \\
\frak{X}(\Pi)_{\rm irr} := \frak{X}(\Pi) \setminus (\iota_{\Pi}^a)^{-1}(V(\Delta(\Pi))) = \frak{X}(\Pi)_{\iota_{\Pi}(\Delta)}.
\end{array}$$
In fact, it is shown ([19, 4.1], [12, $\S 3$]) that $\frak{p}\in \frak{X}(\Pi)$ belongs to $\frak{X}(\Pi)_{\rm irr}$ if and only if $\rho_{\frak{p}}$ is an absolutely irreducible representation. Here a representation $\rho : \Pi \rightarrow {\rm SL}_2(A)$ with a commutative ring $A$ is said to be absolutely irreducible if the composite of $\rho$ with the natural map ${\rm SL}_2(A) \rightarrow {\rm SL}_2(\kappa(\frak{p}))$ is absolutely irreducible over the residue field $\kappa(\frak{p}) = A_{\frak{p}}/\frak{p}A_{\frak{p}}$ for any $\frak{p} \in {\rm Spec}(A)$. 
\\
\\
{\bf Theorem 3.1.1} ([19, 4.3], [12, Corollary 6.8]). {\em The restriction of $\iota_{\Pi}^a$ to $\frak{X}(\Pi)_{\rm irr}$ gives an isomorphism}
$$ \frak{X}(\Pi)_{\rm irr} \simeq \frak{X}^{\rm skein}(\Pi)_{\rm irr}.$$
{\em In terms of algebras, $\iota_{\Pi}$ induces  an isomorphism between $\frak{C}(\Pi)$ and $\frak{B}(\Pi)$  if $\Delta$ is inverted}:
$$ \frak{C}(\Pi)_{\Delta} \simeq  \frak{B}(\Pi)_{\iota_{\Pi}(\Delta)}.$$
\\
{\bf Corollary 3.1.2.} {\it Let  $\rho : \Pi \rightarrow {\rm SL}_2(k)$ be an irreducible representation  and let $[\rho] \in {\cal X}(\Pi)$ also denote the corresponding maximal ideal of $\frak{B}(\Pi)$. Then we have an isomorphism of local rings}:
$$ \frak{C}(\Pi)_{\iota_{\Pi}^a([\rho])} \simeq \frak{B}(\Pi)_{[\rho]}.$$
\\
3.2. {\em The relation between the universal deformation ring and the character scheme.}  Let $k$ be  an algebraically closed field with ${\rm char}(k) \neq 2$ and let ${\cal O}$ be a discrete valuation ring  with residue field $k$. Let $\overline{\rho} : \Pi \rightarrow {\rm SL}_2(k)$ be an irreducible representation and let $\overline{T} : \Pi \rightarrow k$  be a pseudo-${\rm SL}_2$-representation over $k$ given by the character ${\rm tr}(\overline{\rho})$. Let ${\boldsymbol R}_{\overline{T}} (= {\boldsymbol R}_{\overline{\rho}})$ be the universal deformation ring of $\overline{T}$ (or $\overline{\rho}$) as in Sections 1 and 2. Recall that  the universal deformation ring ${\boldsymbol R}_{\overline{T}}$ is a complete local ${\cal O}$-algebra whose  residue field  is  $k$.  On the other hand, let $\frak{B}(\Pi)$ be the universal ${\rm SL}_2$-character algebra of $\Pi$ over $k$.  Then we have the following\\
\\
{\bf Theorem 3.2.1.} {\em Assume that $\overline{\rho}$ is irreducible and let $[\overline{\rho}]$ denote the corresponding maximal ideal of $\frak{B}(\Pi)$. We have an isomorphism of $k$-algebras}
$$ {\boldsymbol R}_{\overline{T}} \otimes_{\cal O} k \; \simeq \; \frak{B}(\Pi)_{[\overline{\rho}]}^{\wedge},$$
{\em where $\frak{B}(\Pi)_{[\overline{\rho}]}^{\wedge}$ denotes the $[\overline{\rho}]$-adic completion of $\frak{B}(\Pi)$.  So, the universal deformation ring can be considered as an infinitesimal deformation of the universal character algebra. For the case that ${\rm char}({\cal O}) = {\rm char}(k)$, we have an isomorphism of ${\cal O}$-algebras
$$ {\boldsymbol R}_{\overline{T}} \; \simeq \; \frak{B}(\Pi)_{[\overline{\rho}]}^{\wedge} \hat{\otimes}_{k} {\cal O},$$
where ${\cal O}$ is considered as a $k$-algebra by the natural inclusion $k \hookrightarrow {\cal O}$.}\\
\\
{\em Proof.} By the construction of ${\boldsymbol R}_{\overline{T}}$ in the proof of Theorem 1.2.1, we have
$$ {\boldsymbol R}_{\overline{T}} = {\cal O}[[ X_g \; (g \in \Pi)]]/{\cal I},$$
where ${\cal I}$ is the ideal of the power series ring ${\cal O}[[ X_g \; (g \in \Pi)]]$ generated by elements of the form: setting $T_g := X_g + \varphi(\overline{T}(g))$, $\varphi$ being the Teichm\"{u}ller lift,
\vspace{.18cm}\\
(1) $T_1 -2$,\\
(2) $T_{g_1g_2} - T_{g_2g_1}$,\\
(3) $T_{g_1}T_{g_2}T_{g_3} + T_{g_1g_2g_3} + T_{g_1g_3g_2} -T_{g_1g_2}T_{g_3} - T_{g_2g_3}T_{g_1} - T_{g_1g_3}T_{g_2}$,\\
(4) $T_g^2 - T_{g^2} -2,$
\vspace{.18cm}\\
where $g, g_1, g_2, g_3 \in \Pi$. 

On the other hand, since the maximal ideal $[\overline{\rho}]$ of $\frak{B}(\Pi)$ corresponds to the maximal ideal $(t_g - \overline{T}(g)\; (g \in \Pi))$ of $\frak{C}(\Pi)$, Corollary 3.1.2 yields
$$ \frak{B}(\Pi)_{[\overline{\rho}]}^{\wedge} \simeq k[[ x_g \; (g \in \Pi)]]/I^{\wedge},$$
where $x_g := t_g - \overline{T}(g) \; (g \in \Pi)$ and $I^{\wedge}$ is the ideal of the power series ring $k[[ x_g \; (g \in \Pi)]]$ generated by  elements of the form $t_1 -2,\;  t_{g_1}t_{g_2} - t_{g_1g_2} - t_{g_1^{-1}g_2}\;  (g_1, g_2 \in \Pi).$ So, in order to show that the correspondence $X_g \otimes 1 \mapsto x_g$ (resp. $X_g  \mapsto x_g \otimes 1$) gives the desired isomorphism ${\boldsymbol R}_{\overline{T}} \otimes_{\cal O} k  \simeq  \frak{B}(\Pi)_{[\overline{\rho}]}^{\wedge}$ (resp. ${\boldsymbol R}_{\overline{T}}   \simeq  \frak{B}(\Pi)_{[\overline{\rho}]}^{\wedge}\hat{\otimes}_k {\cal O}$ for the case that ${\rm char}({\cal O}) = {\rm char}(k)$), it suffices to show the following\\
\\
{\bf Lemma 3.2.2.} {\em Let $T$ be a function on $\Pi$ with values in an integral domain whose characteristic is not $2$. Let ${\rm (P)}$ be the relations given by}
\vspace{.12cm}\\
(P1) $T(1) = 2,$\\
(P2) $T(g_1g_2) = T(g_2g_1)$,\\
(P3) $T(g_1)T(g_2)T(g_3)+T(g_1g_2g_3)+T(g_1g_3g_2)-T(g_1g_2)T(g_3)-T(g_2g_3)T(g_1)-T(g_1g_3)T(g_2) =0 $,\\
(P4) $T(g)^2 - T(g^2) = 2$,
\vspace{.12cm}\\
{\em and let ${\rm (C)}$ be the relations given by}
\vspace{.12cm}\\
(C1) $T(1) = 2$,\\
(C2) $T(g_1)T(g_2) = T(g_1g_2) + T(g_1^{-1}g_2)$,
\vspace{.12cm}\\
{\em where $g, g_1, g_2, g_2, g_3$ are any element in $\Pi$. 

Then ${\rm (P)}$ and ${\rm (C)}$ are equivalent.}\\
\\
{\em Proof of} Lemma 3.2.2.  (P) $\Rightarrow$ (C): Letting $g_2 = g_1$ in (P3), we have
$$T(g_1)^2T(g_3)-T(g_1^2)T(g_3)+T(g_1^2g_3)+T(g_1g_3g_1)-2T(g_1g_3)T(g_1)=0.$$
Using (P2) and (P4), we have
$$2(T(g_3)+T(g_1^2g_3)-T(g_1g_3)T(g_1))=0.$$
Letting $g_3$ be replaced by $g_1^{-1}g_2$ in the above equation and noting $T$ has the value in an integral domain whose characteristic is not $2$, we obtain (C2).

(C) $\Rightarrow$ (P). Letting $g_2=1$ in (C2) and using (C1), we have 
$$ T(g) = T(g^{-1}) \; \mbox{for any}\; g \in \Pi. $$
Exchanging $g_1$ and $g_2$ in (C2) each other and using the above relation, we have
$$ T(g_2)T(g_1) = T(g_2g_1) + T(g_2^{-1}g_1) = T(g_2g_1) + T(g_1^{-1}g_2)$$
and hence we obtain (P2). Next letting $g_1$ be replaced by $g_1g_3$ in (C2), we have
$$ -T(g_1g_3)T(g_2) + T(g_1g_3g_2) + T(g_3^{-1}g_1^{-1}g_2) = 0, \leqno{(3.2.2.1)} $$
and letting $g_2$ be replaced by $g_2g_3$ in (C2), we have
$$ -T(g_1)T(g_2g_3) + T(g_1g_2g_3) + T(g_1^{-1}g_2g_3) = 0. \leqno{(3.2.2.2)}$$
By (C2), we have
$$ \begin{array}{ll}
T(g_3^{-1}g_1^{-1}g_2) & = T(g_3)T(g_1^{-1}g_2) - T(g_3g_1^{-1}g_2) \\
 & = T(g_3)T(g_1)T(g_2)-T(g_1g_2)T(g_3) - T(g_3g_1^{-1}g_2).
\end{array}$$
Hence, using (P2) proved already, we have
$$ \begin{array}{ll}
T(g_3^{-1}g_1^{-1}g_2) + T(g_1^{-1}g_2g_3) & = T(g_1)T(g_2)T(g_3)-T(g_1g_2)T(g_3) \\
                                                                           & \;\;\;\;\;\;\;\; -T(g_3g_1^{-1}g_2)+T(g_1^{-1}g_2g_3) \\
& = T(g_1)T(g_2)T(g_3)-T(g_1g_2)T(g_3). 
\end{array} \leqno{(3.2.2.3)} $$
Summing up (3.2.2.1) and (3.2.2.2) together with (3.2.2.3), we obtain (P3). Finally putting $g_1 = g_2$ in (C2) and using (C1), we obtain (P4). $\;\;\Box$\\

By Lemma 2.2.3 and Theorem 3.2.1, we have the following\\
\\
{\bf Corollary 3.2.3.} {\em Assume  that $[\overline{\rho}]$ is a regular point of the scheme $\frak{X}(\Pi)$, namely, $\frak{B}(\Pi)_{[\overline{\rho}]}$ is a regular local ring. Then the dimension $d$ of the relative tangent space $\frak{t}_{{\boldsymbol R}_{\overline{T}}/{\cal O}}$ of ${\boldsymbol R}_{\overline{T}}$ is equal to  the dimension of the irreducible componennt of $\frak{X}(\Pi)$ containing $[\overline{\rho}]$, and $\frak{B}(\Pi)_{[\overline{\rho}]}^{\wedge}$ is a power series ring over $k$ on a regular system of parameters  $x_1,\dots,x_d$. Hence we have  a surjective ${\cal O}$-algebra homomorphism 
$$ \lambda : {\cal O}[[X_1,\dots, X_d]] \longrightarrow {\boldsymbol R}_{\overline{T}}$$
 in ${\cal C}$ such that the image of $\lambda(X_i)$ in ${\boldsymbol R}_{\overline{T}}\otimes_{\cal O} k \simeq \frak{B}(\Pi)_{[\overline{\rho}]}^{\wedge}$ is $x_i$ $(1\leq i \leq d)$.}\\
\\

\begin{center}
{\bf 4. Examples for 2-bridge knot groups}
\end{center}

In this section, we investigate examples concerning Riley representations of 2-bridge knot groups. 

In Subsection 4.1, we recall some results on the Riley representations of 2-bridge knot groups. We refer to [1] for basic information on 2-bridge knots and  [17], [18] for the details on Riley representations.
In Subsection 4.2, we describe the character scheme/variety of ${\rm SL}_2$-representations of a 2-bridge knot group. For this, we refer to [6]. In Subsection 4.3, we give an explicit form of the universal  deformation of a Riley representation. \\

4.1.  {\em $2$-bridge knots and Riley representations.} Let $K$ be a 2-bridge knot in the 3-sphere $S^3$, given as the Schubert form $\frak{b}(m,n)$ where $m$ and $n$ are odd integers with $m>0, -m<n<m$ and ${\rm g.c.d}(m,n) =1$. Let $\Pi_K$ be the knot group $\pi_1(S^3 \setminus K)$. The group $\Pi_K$ is known to have a presentation of the form
$$ \Pi_K = \langle a, b \; | \; wa = bw \rangle,$$
where $w$ is a word $w(a,b)$ of $a$ and $b$ which has the following symmetric form
$$ \begin{array}{l}
w = w(a,b) = a^{\epsilon_1} b^{\epsilon_2} \cdots a^{\epsilon_{m-2}} b^{\epsilon_{m-1}}, \\
\epsilon_i = (-1)^{[in/m]} = \epsilon_{m-i} \;\; ([\, \cdot \,] = \,\mbox{Gauss symbol}). \end{array}$$
Let $F$ be the free group generated by $a$ and $b$, and let $\pi : F \rightarrow \Pi_K$ be the natural homomorphism.

Let $A$ be  a commutative ring with identity. For $\alpha \in A^{\times}$ and $\beta \in A$, we consider two matrices $C(\alpha)$ and $D(\alpha,\beta)$ in ${\rm SL}_2(A)$ defined by
$$ C(\alpha) := \left(  \begin{array}{cc} \alpha & 1 \\ 0 & \alpha^{-1} \end{array} \right), \;\;  D(\alpha,\beta) := \left(  \begin{array}{cc} \alpha & 0 \\ \beta & \alpha^{-1} \end{array} \right)   \leqno{(4.1.1)}$$
and we set
$$ W(\alpha,\beta) := C(\alpha)^{\epsilon_1} D(\alpha,\beta)^{\epsilon_2}\cdots C(\alpha)^{\epsilon_{m-2}} D(\alpha,\beta)^{\epsilon_{m-1}}.$$
It is easy to see that there are (Laurent) polynomials $w_{ij}(t,u) \in \mathbb{Z}[t^{\pm},u]$ $(1\leq i,j \leq 2)$ such that $W(\alpha,\beta) = (w_{ij}(\alpha,\beta))$. Let $f_{(\alpha,\beta)} : F \rightarrow {\rm SL}_2(A)$
be the homomorphism defined by
$$ f_{(\alpha,\beta)}(a) := C(\alpha), \;\; f_{(\alpha,\beta)}(b) := D(\alpha,\beta).$$
We call a representation
$$ \rho : \Pi_K \longrightarrow {\rm SL}_2(A)$$
the {\em Riley representation over} $A$ {\em of type} $(\alpha,\beta)$, denoted by $\frak{r}_{(\alpha,\beta)}$, if $f_{(\alpha,\beta)}$ factors through $\rho$, namely, $\rho \circ \pi = f_{(\alpha,\beta)}$.

Let $k$ be an algebraically closed field. The following Theorem 4.1.2 was proved by Riley [17], [18] for the case where $k$ is the field of complex numbers. The proof therein works as well for any algebraically closed field. \\
\\
{\bf Theorem 4.1.2} ([17], [18]).  {Let $\varphi(t,u) := w_{11}(t,u) + (t^{-1} - t)w_{12}(t,u) \in \mathbb{Z}[t^{\pm},u]$.}\\
(1)  {\em There is a unique polynomial $\Phi(x,u) \in \mathbb{Z}[x,u]$ such that}
$$ \Phi(t+t^{-1}, u) = t^{l} \varphi(t,u)$$
{\em for an integer $l$.}\\
(2) {\em The homomorphism $f_{(\alpha,\beta)}$ $(\alpha \in k^{\times}, \beta \in k)$ factors through the Riley representation $\frak{r}_{(\alpha,\beta)}$ over $k$ if and only if we have}
$$ \Phi(\alpha + \alpha^{-1}, \beta) = 0.$$
(3) {\em Any non-Abelian ${\rm SL}_2(k)$-representation of $\Pi_K$ is equivalent to a Riley representation $\frak{r}_{(\alpha,\beta)}$ for some $\alpha \in k^{\times}$ and $\beta \in k$.}
\\
\\
For the properties of the polynomial $\Phi(x,u)$,  Riley showed, among others, the following\\
\\
{\bf Proposition 4.1.3} ([17]).  {\em The polynomial $\Phi(2,u) = \varphi(1,u) \in \mathbb{Z}[u]$ is monic up to multiplication by $\pm 1$ and its discriminant ${\rm disc}(\Phi(2,u))$ is an odd integer. If ${\rm char}(k)$ does not divide ${\rm disc}(\Phi(2,u))$, then any root of $\Phi(2,u) = 0$ in $k$ is non-zero and simple. }\\
\\
By Hensel's lemma, we have the following\\
\\
{\bf Corollary 4.1.4.} {\em Let ${\cal O}$ be a complete discrete valuation ring with residue field $k$. For any root $\beta$ of $\Phi(2,u) = 0$ in $k$, there is a unique power series $u(x) \in {\cal O}[[x-2]]^{\times}$ such that $\beta = u(2) \, {\rm mod}\, \frak{m}_{\cal O}$ and $\Phi(x,u(x)) \equiv 0$ in ${\cal O}[[x-2]]$.}\\
\\
{\bf Example 4.1.5.}  (1) Let $K$ be the trefoil $\frak{b}(3,1)$. Then we have $w = ab$ and $\varphi(t,u) = t^2(t^2+t^{-2}+u-1)$. Hence  $\Phi(x,u) = x^2  + u -3$ and $\Phi(2,u) = u+1$. Therefore $\beta = -1$ and $u(x) = 3 - x^2$ for any $k$.\\
(2) Let $K$ be the figure eight $\frak{b}(5,3)$. Then we have $w = ab^{-1}a^{-1}b$ and $\varphi(t,u) = u^2 + (t^2 + t^{-2} -3)u - (t^2 + t^{-2} -3)$. Hence $\Phi(x,u) = u^2 +(x^2 - 5)u -(x^2 - 5)$ and
$\Phi(2,u) = u^2 -  u + 1$. Therefore, if ${\rm char}(k) \neq 2, 3$, then $\beta = \displaystyle{ \frac{1}{2}(1 \pm \sqrt{-3})} \in k$ and $u(x) = \displaystyle{\frac{1}{2} \{5 -x^2 \pm \sqrt{(x^2 -1)(x^2 - 5)} \}} \in {\cal O}[[x-2]]^{\times}$ where $\sqrt{(x^2 -1)(x^2 - 5)}$ stands for an element of ${\cal O}[[x-2]]$ whose square is $(x^2 -1)(x^2 - 5)$\\

4.2. {\em Character varieties}.  We keep the notations  in Subsection 3.1. Let $k$ denote an algebraically closed field and let ${\cal X}(\Pi_K)$ denote the ${\rm SL}_2(k)$-character variety of $\Pi_K$. The proof of Proposition 1.4.1 of [3] tells us that any $\tau_g$ ($g \in \Pi_K$) is given as a polynomial of $\tau_a (= \tau_b)$ and $\tau_{ab}$ with coefficients in $\mathbb{Z}$. In particular, the coordinate ring $k[{\cal X}(\Pi_K)]$ is generated by $\tau_a$ and $\tau_{ab}$. We let $x$ and $y$ denote the variables corresponding,  respectively,  to the coordinate functions  $\tau_a$ and $\tau_{ab}$ on ${\cal X}(\Pi_K)$ embedded in $k^2$. This variable $x$ is consistent with the variable $x$ of $\Phi(x,u)$ in Theorem 4.1.2 (and so causes no confusion). In fact, the coordinate variables $x$ and $y$ are related with $t$ and $u$ in Theorem 4.1.2 by 
$$x = t + t^{-1}, \;\; y = t^2 + t^{-2} + u = x^2 + u -2,$$
since we have
$$\tau_a(\frak{r}_{(\alpha,\beta)}) = {\rm tr}(C(\alpha)) = \alpha + \alpha^{-1}, \;\; \tau_{ab}(\frak{r}_{(\alpha,\beta)}) = {\rm tr}(C(\alpha)D(\alpha,\beta)) = \alpha^2 + \alpha^{-2} + \beta.$$

By Theorem 4.1.2 (2), (3), characters of irreducible ${\rm SL}_2(k)$-representations of $\Pi_K$ correspond bijectively to points on the algebraic curve in $k^2$ defined by the equation
$$ \Phi(x, y - x^2 +2) = 0,$$
except the finitely many intersection points with the algebraic curve $y -x^2 +2 = 0$. The points on the latter curve $y -x^2 +2 = 0$ correspond (not bijectively) to characters of reducible ${\rm SL}_2(k)$-representations of $\Pi_K$. 
 It is also shown ([6, Proposition 3.4.1]) that the ideal generated by $\Phi(x, y - x^2 +2)$ in $k[x,y]$ is a radical ideal. Thus we have the following\\
\\
{\bf Theorem 4.2.1} ([6, Theorem 3.3.1]). {\em The character variety ${\cal X}(\Pi_K)$ is the  affine algebraic curve in $k^2$ defined by the equation}
$$ (y -x^2 +2)\Phi(x, y-x^2+2) = 0,$$
and the coordinate ring  of ${\cal X}(\Pi_K)$ is given by 
$$ k[{\cal X}(\Pi_K)] \simeq k[x,y]/((y - x^2 + 2)\Phi(x, y-x^2+2)).$$
{\em Here the points on the algebraic curve $\Phi(x, y-x^2+2) = 0$ correspond bijectively to irreducible ${\rm SL}_2(k)$-characters of $\Pi_K$ except the finitely many intersection points with the algebraic curve $y -x^2 +2=0$, and the points on 
$y -x^2 +2=0$ correspond $($not bijectively$)$ to reducible ${\rm SL}_2(k)$-characters of $\Pi_K$.}\\ 
\\ 
 {\bf Example 4.2.2.} (1) When $K$ is the trefoil $\frak{b}(3,1)$, we see $\Phi(x, y-x^2+2) = y-1$. Hence ${\cal X}(\Pi_K)$ is given by $(y - x^2+2)(y-1)=0$.\\
 (2) When $K$ is the figure eight $\frak{b}(5,3)$, we have $\Phi(x,y-x^2+2) = y^2 -(1+x^2)y + 2x^2 -1$. Hence ${\cal X}(\Pi_K)$ is given by $(y- x^2 + 2)\{ y^2 -(1+x^2)y + 2x^2 -1   \} = 0$. \\
 \\
 Przytycki and Sikora proved the following theorem for the case where $k$ is the field of complex numbers. Their proof works well for any algebraically closed field whose characteristic is not $2$. 
 \\
 \\
 {\bf Theorem 4.2.3} ([16, Theorem 7.3]). {\em Assume ${\rm char}(k) \neq 2$. Then the universal ${\rm SL}_2$-character algebra $\frak{B}(\Pi_K)$ is reduced and hence $ \frak{B}(\Pi_K) = k[{\cal X}(\Pi_K)].$}\\
 \\ 
4.3.  {\em The universal deformation}.  As in Subsection 3.2, let $k$ be an algebraically closed field with ${\rm char}(k) \neq 2$ and  let ${\cal O}$ be  a complete discrete valuation ring  with residue field $k$. We assume further that ${\rm char}(k)$ does not divide the discriminant of $\Phi(2,u) \in \mathbb{Z}[u]$. 

Let $\overline{\rho} : \Pi_K \rightarrow {\rm SL}_2(k)$ be a Riley representation $\frak{r}_{(1,\beta)}$ so that
$$ \overline{\rho}(a) := C(1), \;\; \overline{\rho}(b) := D(1,\beta),$$
where $\beta$ is a root of $\Phi(2,\beta) = 0$. By Proposition 4.1.3, $\overline{\rho}$ is  irreducible and $(x,y) = (2,\beta+2)$ is a non-singular point on ${\cal X}(\Pi_K)$. 

Let $u(x)$ be the power series in  Corollary  4.1.4 and set $\overline{u}(x) := u(x) \, {\rm mod}\, \frak{m}_{\cal O}$. By Theorem 4.2.1, we have the isomorphism
$$ k[{\cal X}(\Pi_K)]_{[\overline{\rho}]}^{\wedge}\simeq (k[x,y]/(\Phi(x,y-x^2+2)))_{(x-2,y-(\beta+2))}^{\wedge} \simeq k[[x-2]], \leqno{(4.3.1)}$$
where the second isomorphism is given by $y \mapsto x^2 + \overline{u}(x) -2$. So $x-2$ is a local parameter of ${\cal X}(\Pi_K)$ at $[\overline{\rho}]$.

Let $({\boldsymbol R}_{\overline{\rho}}, {\boldsymbol \rho})$ be the universal deformation of $\overline{\rho}$. By Theorem 3.2.1, Theorem 4.2.3 and (4.3.1), we have 
$${\boldsymbol R}_{\overline{\rho}}\otimes_{\cal O} k \simeq k[[x-2]],$$
where $T_a \, \mbox{mod} \, {\cal I} = {\rm tr}({\boldsymbol \rho}(a)) \in {\boldsymbol R}_{\overline{\rho}}$ corresponds to $x$.  By Corollary 3.2.3, we have \\
\\
{\bf Lemma 4.3.2.} {\em The dimension of the relative tangent space $\frak{t}_{{\boldsymbol R}_{\overline{\rho}}/{\cal O}}$ is $1$ and there is a surjective ${\cal O}$-algebra homomorphism 
$$\lambda : {\cal O}[[X]] \longrightarrow {\boldsymbol R}_{\overline{\rho}}$$
in ${\cal C}$ such that $\lambda(X) = {\rm tr}({\boldsymbol \rho}(a)) - 2$.}\\
\\
 In the following Theorem 4.3.3, we show that the map $\lambda$ in Lemma 4.3.2 is in fact an isomorphism, and we give an explicit form of the universal deformation $({\boldsymbol R}_{\overline{\rho}}, {\boldsymbol \rho})$. We remark on the notation used in the following: For $p(x) \in {\cal O}[[x-2]]$ with $p(2) \in {\cal O}^{\times}$, $\sqrt{p(x)}$ stands for an element in ${\cal O}[[x-2]]$ whose square is $p(x)$. If $p(2) = 1$, we adopt the unique one normalized by $\sqrt{p(2)} = 1$. For $p(x) \in {\cal O}[[x-2]]$ with $p(2) \in \frak{m}_{\cal O}$, $\sqrt{p(x)}$ is an element of a quadratic extension of $E((x-2))$ whose square is  $p(x)$, where $E$ is the field of fractions of ${\cal O}$.
In particular, $\sqrt{x-2}$ is a prime element of a quadratic extension $L$ of $E((x-2))$ and we denote by ${\cal O}[[\sqrt{x-2}]]$ the integral closure of ${\cal O}[[x-2]]$ in $L$. For $p(x) \in {\cal O}[[x-2]]$ with $p(2)=0$, 
  we may write  $p(x) = (x-2)p_1(x)$ with $p_1(x) \in {\cal O}[[x-2]]$ and so we have $\frac{p(x)}{\sqrt{x-2}} = \sqrt{x-2}p_1(x) \in {\cal O}[[\sqrt{x-2}]]$. 
\\
\\
{\bf Theorem 4.3.3.} {\em We let $v(x) := \sqrt{1+\frac{x^2-4}{u(x)}} \in {\cal O}[[x-2]]^{\times}$ and define $U(x) \in {\rm SL}_2({\cal O}[[\sqrt{x-2}]])$ by}
$$ U(x) := \left(   
\begin{array}{cc}
\displaystyle{\frac{1}{\sqrt{v(x)}}} & \displaystyle{\frac{1-v(x)}{\sqrt{v(x)}\sqrt{x^2-4}}}\\
\displaystyle{\frac{\sqrt{x^2-4}}{2\sqrt{v(x)}}} & \displaystyle{\frac{1+v(x)}{2\sqrt{v(x)}}}
\end{array}
\right).
$$
{\em We  define $A(x), B(x) \in {\rm SL}_2({\cal O}[[x-2]])$ by}
$$ \begin{array}{l} A(x) := U(x)C(t)U(x)^{-1} = \left( \begin{array}{cc} \displaystyle{\frac{x}{2}} & 1 \\ \displaystyle{\frac{x^2-4}{4}} & \displaystyle{\frac{x}{2}} \end{array}  \right), \\
B(x) := U(x)D(t,u(x))U(x)^{-1} = \left( \begin{array}{cc} \displaystyle{\frac{x}{2}} & \displaystyle{\frac{(1-v(x))^2u(x)}{x^2-4}} \\  \displaystyle{\frac{(1+v(x))^2u(x)}{4}} & \displaystyle{\frac{x}{2}} \end{array}  \right),
\end{array}
$$
{\em where $t$ is an  element ${\cal O}[[\sqrt{x-2}]]$ such that $t + t^{-1} = x$ and $C(t), D(t,u(x))$ are the matrices over ${\cal O}[[\sqrt{x-2}]]$ defined in $(4.1.1)$. We define the deformation of $\overline{\rho}$}
$$ {\boldsymbol \rho^u} : \Pi_K \longrightarrow {\rm SL}_2({\cal O}[[x-2]])$$
{\em by}
$$ {\boldsymbol \rho^u}(a) := A(x), \;\; {\boldsymbol \rho^u}(b) := B(x).$$
{\em Then there is an isomorphism $\psi : {\boldsymbol R}_{\overline{\rho}} \stackrel{\sim}{\rightarrow} {\cal O}[[x-2]]$ in ${\cal C}$ such that $\psi \circ {\boldsymbol \rho} \approx {\boldsymbol \rho^u}$.}
\\
\\
{\em Proof.}  Firstly, let us check that $U(x) \in {\rm SL}_2({\cal O}[[\sqrt{x-2}]])$ and $A(x), B(x) \in {\rm SL}_2({\cal O}[[x-2]])$. Since $v(2) = 1$, $\sqrt{v(x)} \in {\cal O}[[x-2]]^{\times}$ and $1-v(x) = (x-2)p(x)$ with some $p(x) \in {\cal O}[[x-2]]$ and hence all entries of $U(x)$ are lying in ${\cal O}[[\sqrt{x-2}]]$. Further we easily see  $\det U(x) =1$ and also $U(2) = I$. As for $A(x), B(x)$, we see immediately that $\det A(x) = \det C(t) = 1$ and $\det B(x) = \det D(t,u(x)) = 1$. The straightforward computations of $U(x)C(t)U(x)^{-1}$ and $U(x)D(t,u(x))U(x)^{-1}$  using $t+ t^{-1} = x, t - t^{-1} = \sqrt{x^2 -4}, x^2 -4 + u(x) = v(x)^2 u(x)$ yield the desired matrices in the statement, from which  we easily see that all entries of $A(x)$ and $B(x)$ are  lying in ${\cal O}[[x-2]]$.

Next, let us show that ${\boldsymbol \rho^u}$ is a deformation of $\overline{\rho}$ over ${\cal O}[[x-2]]$. Since ${\boldsymbol \rho^u}$ is equivalent to the Riley representation $\frak{r}_{(t,u(x))}$ over ${\cal O}[[\sqrt{x-2}]]$, ${\boldsymbol \rho^u}$ is indeed a representation. Since $A(2) \, \mbox{mod}\, \frak{m}_{{\cal O}[[x-2]]} = C(1)$ and $B(2) \, \mbox{mod}\, \frak{m}_{{\cal O}[[x-2]]} = D(1,\beta)$, we find that ${\boldsymbol \rho^u}$ mod $\frak{m}_{{\cal O}[[x-2]]} = \overline{\rho}$.

Finally, by the universality of $({\boldsymbol R}_{\overline{\rho}}, {\boldsymbol \rho})$, we have a homomorphism $\psi : {\boldsymbol R}_{\overline{\rho}}  \rightarrow {\cal O}[[x-2]]$ in ${\cal C}$ such that $\psi \circ {\boldsymbol \rho} \approx {\boldsymbol \rho^u}$. So we have $\psi({\rm tr}({\boldsymbol \rho}(a))) = {\rm tr}({\boldsymbol \rho^u}(a)) = x$. On the other hand, by Lemma 4.3.2, we have $\lambda(X) = {\rm tr}({\boldsymbol \rho}(a)) - 2$. Therefore $\psi \circ \lambda (X) = x - 2$ and hence $\psi \circ \lambda$ is an isomorphism ${\cal O}[[X]] \simeq {\cal O}[[x-2]]$. Since $\lambda$ is surjective, $\psi$ and $\lambda$ must be isomorphisms in ${\cal C}$.  $\;\; \Box$
\\
\\
{\bf Remark 4.3.4.} (1) By the construction of ${\boldsymbol \rho^u}$, if $\overline{\rho}$ is defined over a subfield $k'$ of $k$,  the representation ${\boldsymbol \rho^u}$ is also defined over a ring ${\cal O}'[[x-2]]$, where ${\cal O}'$ is a complete discrete valuation ring with residue field $k'$. For example, if $\overline{\rho}$ is defined over a prime field $\mathbb{F}_p$ of $p$ elements, ${\boldsymbol \rho^u}$ can be defined over $\mathbb{Z}_p[[x-2]]$, where $\mathbb{Z}_p$ is the ring of $p$-adic integers.\\
(2) Suppose  ${\rm char}({\cal O}) = {\rm char}(k)$ so that  ${\cal O} = k[[\hbar]]$. Then the representation ${\boldsymbol \rho^u}$ is independent of $\hbar$. However, there are deformations of $\overline{\rho}$ which depend on $\hbar$. For example,  letting $\rho_n(a) := A((1+\hbar)^n + (1+\hbar)^{-n})$ and $\rho_n(b) := B((1+\hbar)^n + (1+\hbar)^{-n})$ for $n \in \mathbb{Z}$, we have a family of deformations $\rho_n$ of $\bar{\rho}$ over $k[[\hbar]]$. \\
\\
{\bf Example 4.3.5.} (1) Let $K$ be the trefoil $\frak{b}(3,1)$ and assume ${\rm char}(k) \neq 2$. We then have $\beta = -1, u(x) = 3-x^2$  and $v(x) = \displaystyle{\frac{1}{\sqrt{x^2 -3}}}$. For example, we can consider $\overline{\rho} = \frak{r}_{(1,-1)} : \Pi_K \rightarrow {\rm SL}_2(\mathbb{F}_p)$ for an odd prime number $p$, where $\mathbb{F}_p \subset k$. Then ${\boldsymbol \rho^u}$ can be a representation into ${\rm SL}_2(\mathbb{Z}_p[[x-2]])$, strictly equivalent to ${\boldsymbol \rho}$ over ${\cal O}[[x-2]]$. \\
(2) Let $K$ be the figure eight $\frak{b}(5,3)$ and assume ${\rm char}(k) \neq 2, 3$. We then have $\beta = \frac{1}{2}(1 \pm \sqrt{-3}), u(x) = \frac{1}{2} \{5 -x^2 \pm \sqrt{(x^2 -1)(x^2 - 5)} \}$ and $v(x)^2 = \frac{1}{2}\{ x^2 -2 \pm \frac{(x^2-4)\sqrt{x^2-1}}{\sqrt{x^2-5}}\}$. For example, we can consider $\overline{\rho} = \frak{r}_{(1,\beta)} : \Pi_K \rightarrow {\rm SL}_2(\mathbb{F}_p(\beta))$ for $p \neq 2, 3$, where $\mathbb{F}_p(\beta) \subset k$. Then ${\boldsymbol \rho^u}$ can be a representation into ${\rm SL}_2(\mathbb{Z}_p[\beta][[x-2]])$, strictly equivalent to ${\boldsymbol \rho}$ over ${\cal O}[[x-2]]$. 
\\ 
\\
\begin{center}
{\bf 5. The universal deformation of a holonomy representation}
\end{center}

In this section, we apply our deformation theory to the case where $\Pi$ is the fundamental group of the complement of a hyperbolic knot in the 3-sphere and $\overline{\rho}$ is (a lift of) the holonomy representation. 

In Subsection 5.1, we recall Thurston's theorem on deformations of hyperbolic structures ([21]). In Subsection 5.2, we then describe the universal deformation of $\overline{\rho}$ by using Thurston's theorem, and discuss some analogies with $p$-adic Galois deformations.

In this section, we work over the field $k = \mathbb{C}$ of complex numbers.\\

5.1. {\em Holonomy representation and Thurston's theorem.}  Let $K$ be a hyperbolic knot in the $3$-sphere $S^3$ and let $\Pi_K := \pi_1(S^3 \setminus K)$ be the knot group. The complement $S^3 \setminus K$ is a complete hyperbolic  $3$-manifold of finite volume with a cusp, which is given as a quotient of the hyperbolic $3$-space $\mathbb{H}^3$ by a discrete, torsion free subgroup of ${\rm PSL}_2(\mathbb{C}) =  {\rm SL}_2(\mathbb{C})/\{ \pm I \} = \mbox{Aut}(\mathbb{H}^3)$. To the complete hyperbolic structure on $S^3 \setminus K$ we can associate a faithful representation $\rho_{\rm hol} : \Pi_K \rightarrow {\rm PSL}_2(\mathbb{C})$, called the {\em holonomy representation}. The holonomy representation $\rho_{\rm hol}$ can be lifted to an ${\rm SL}_2(\mathbb{C})$-representation, and thus we fix such a lift
$$\overline{\rho}_{\rm hol} : \Pi_K \longrightarrow {\rm SL}_2(\mathbb{C}),$$
which is known to be irreducible.

Let ${\cal X}(\Pi_K)$ be the ${\rm SL}_2(\mathbb{C})$-character variety of $\Pi_K$ as defined in Subsection 3.1 and let ${\cal X}(\Pi_K)^{\rm hol}$ be the irreducible component of ${\cal X}(\Pi_K)$ containing $[\overline{\rho}_{\rm hol}]$. We choose any meridian $\mu$ of the knot $K$ and consider the map $\tau_{\mu} : \frak{X}(\Pi_K)^{\rm hol} \rightarrow \mathbb{C}$ defined by $\tau_{\mu}([\rho]) = {\rm tr}(\rho)(\mu)$. Then Thurston has proved the following\\
\\
{\bf Theorem 5.1.1} ([21]). {\em The map $\tau_{\mu}$ is bianalytic in a neighborhood of $[\overline{\rho}_{\rm hol}]$.}\\
\\
Therem 5.1.1 implies that ${\cal X}(\Pi_K)^{\rm hol}$  is a complex algebraic curve and $\tau_{\mu}$ gives a local parameter around the smooth point $[\overline{\rho}_{\rm hol}]$. Hence we have the following\\
\\
{\bf Corollary 5.1.2.} {\em We have the following isomorphism of $\mathbb{C}$-algebras}
$$ \mathbb{C}[{\cal X}(\Pi_K)]_{[\overline{\rho}_{\rm hol}]}^{\wedge} \; \simeq \; \mathbb{C}[[z]],$$
{\em where $z$ is a variable corresponding to $\tau_{\mu} - \tau_{\mu}(\overline{\rho}_{\rm hol}) = \tau_{\mu} -2$.}\\

5.2.  {\em The universal deformation of the holonomy representation.} Let 
$$ \mbox{\boldmath $\rho$}_{\rm hol} : \Pi_K \longrightarrow {\rm SL}_2({\boldsymbol R}_{\overline{\rho}_{\rm hol}})$$
be the universal deformation of $\overline{\rho}_{\rm hol}$ in Theorem 2.2.2, where the universal deformation ring ${\boldsymbol R}_{\overline{\rho}_{\rm hol}}$ is a complete local algebra over ${\cal O} = \mathbb{C}[[\hbar]]$. We assume that the universal ${\rm SL}_2$-character algebra $\frak{B}(\Pi_K)$ is reduced so that $\frak{B}(\Pi_K) = k[{\cal X}(\Pi_K)]$. Then,  by Theorem 3.2.1 and Corollary 5.1.2, we have the following  \\
\\
{\bf Theorem 5.2.1.} {\em Under the above assumption, we have the following isomorphism of ${\cal O}$-algebras}
$${\boldsymbol R}_{\overline{\rho}_{\rm hol}} \; \simeq \; {\cal O}[[z]].$$
\\
{\bf Remark 5.2.2.}  The analogy between the structures of ${\cal X}(\Pi_K)^{\rm hol}$ and the deformation space of nearly ordinary $p$-adic Galois representations was firstly pointed out by Kazuhiro Fujiwara (cf. [11, Chapter 14]). Theorem 5.2.1 may make this analogy more precise as follows: Let $\frak{D}$ be the irreducible component of the deformation space of hyperbolic structures on $S^3 \setminus K$ containing the complete hyperbolic structure, say $z^{\circ}$ . By Thurston's theory on hyperbolic Dehn filling ([21]), a neighborhood of $z^{\circ}$ in $\frak{D}$ is homeomorphic to  a neighborhood of $[\overline{\rho}_{\rm hol}]$ in $\frak{X}(\Pi_K)^{\rm hol}$, associating to an (incomplete) hyperbolic structure  the holonomy representation. So Theorem 5.2.1 gives the isomorphism between the universal deformation ring ${\boldsymbol R}_{\overline{\rho}_{\rm hol}}$ and the complete local ring of $\frak{D}$ at $z^{\circ}$, where the parameter $z$ in Theorem 5.2.1 may also be considered as hyperbolic structure (Dehn filling coefficient). Noting that the restriction of $\rho$ to the peripheral group of $K$ (the fundamental group of the boundary of a tubular neighborhood of $K$) is equivalent to an uppertriangular representation, this isomorphism is quite analogous to the isomorphism between the universal deformation ring for $p$-adic ordinary Galois representatios and the $p$-adic ordinary Hecke algebra, which implies that any $p$-adic ordinary deformation of a given modular Galois representation over $\mathbb{F}_p$ is associated to a $p$-adic ordinary modular form (cf. [5]). Here we may observe that hyperbolic structures (Dehn filling coefficients) correspond to $p$-adic ordinary modular forms ($p$-adic weights).   See also Ohtani's article [14] for a related analogy.\\
\\
\begin{flushleft}
{\small
{\bf References}\\
{[1]} G. Burde,  H. Zieschang, Knots,  Second edition. de Gruyter Studies in Mathematics, {\bf 5}. Walter de Gruyter \& Co., Berlin, 2003.\\
{[2]} H. Carayol, Formes modulaires et repr\'{e}sentations galoisiennes \`{a} valeurs dans un anneau local complet, $p$-adic monodromy and the Birch and Swinnerton-Dyer conjecture (Boston, MA, 1991), 213--237, Contemp. Math., 165, Amer. Math. Soc.,  Providence, RI,  1994. \\
{[3]} M. Culler, P. Shalen, Varieties of group representations and
   splittings of $3$-manifolds, Ann. of Math. {\bf 117}  (1983), 109--146.\\
{[4]} H. Hida, Galois representations into ${\rm GL}_2(\mathbb{Z}_{p}[[X]])$ attached to ordinary cusp forms, Invent. Math. {\bf 85} (1986), 545--613.\\
{[5]} H. Hida, Modular forms and Galois cohomology, Cambridge Studies in Advanced Mathematics, {\bf 69}, Cambridge University Press, Cambridge, 2000. \\
{[6]} T. Le, Varieties of representations and their subvarieties of cohomology jumps for certain knot
groups,  Russian Acad. Sci. Sb. Math. {\bf 78} (1994), 187--209.\\
{[7]} A. Lubotzky, A. Magid, Varieties of representations of finitely generated groups, Mem. Amer. Math. Soc. 58 (1985),  no. 336, xi+117 pp. \\
{[8]} B. Mazur, Remarks on the Alexander polynomial, available at http://www.math.harvard.edu/~mazur/older.html.\\
{[9]} B. Mazur, Deforming Galois representations, Galois groups over $\mathbb{Q}$ (Berkeley, CA, 1987),  385--437,  Math. Sci. Res. Inst. Publ., 16, Springer, New York, 1989.\\ 
{[10]} B. Mazur, The theme of $p$-adic variation, Mathematics: frontiers and
perspectives, 433--459, Amer. Math. Soc., Providence,  RI, 2000. \\
{[11]} M. Morishita, Knots and Primes -- An Introduction to Arithmetic Topology, Universitext. Springer, London, 2012. \\
{[12]} K. Nakamoto, Representation varieties and character varieties, Publ. Res. Inst. Math. Sci.  {\bf 36}  (2000),  no. 2, 159--189.\\
{[13]} L. Nyssen, Pseudo-repr\'{e}sentations, Math. Ann. {\bf 306} (1996), no.2, 257--283.\\
{[14]} S. Ohtani, An analogy between representations of knot groups and Galois groups, Osaka J. Math. {\bf 48}  (2011),  no. 4, 857--872.\\
{[15]} C. Procesi, The invariant theory of $n \times n$ matrices,  Advances in Math.  {\bf 19} (1976) no. 3, 306--381.\\
{[16]} J. Przytycki, A. Sikora, On skein algebras and $Sl_2(\mathbb{C})$-character varieties, Topology {\bf 39}  (2000), 115--148.\\
{[17]} R. Riley, Nonabelian representations of 2-bridge knot groups, Quar. J. Oxford {\bf 35}  (1984), 191--208.\\
{[18]} R. Riley, Holomorphically parameterized families of subgroups of ${\rm SL}(2,\mathbb{C})$, Mathematika {\bf 32} (1985),  no. 2,  248--264.\\
{[19]} K. Saito, Character variety of representations of a finitely generated group in ${\rm SL}_2$, Topology and Teichm\"{u}ller spaces (Katinkulta, 1995),  253--264, World Sci. Publ., River Edge, NJ,  1996.\\
{[20]} R. Taylor, Galois representations associated to Siegel modular forms of low weight, Duke Math. J. {\bf 63} (1991), no. 2, 281--332. \\
{[21]} W. Thurston, The geometry and topology of 3-manifolds, Lect. Note, Princeton, 1979. \\
{[22]} J. Tilouine, Deformations of Galois representations and Hecke algebras, Published for The Mehta Research Institute of Mathematics and 
Mathematical Physics, Allahabad; by Narosa Publishing House, New Delhi, 1996. \\
{[23]} A. Wiles, On ordinary $\lambda$-adic representations associated to modular forms, Invent. Math. {\bf 94} (1988), no.3, 529--573.\\
}
\end{flushleft}

\vspace{2.5cm}

{\small 
Masanori Morishita\\
Faculty of Mathematics, Kyushu University, 744, Motooka, Nishi-ku, Fukuoka, 819-0395, Japan\\
morisita@math.kyushu-u.ac.jp
\vspace{.2cm}\\
Yu Takakura\\
3-3-1-602 Nagaoka, Minami-ku, Fukuoka, 815-0075, Japan\\
ma205017@yahoo.co.jp
\vspace{.2cm}\\
Yuji Terashima\\
Department of Mathematical and Computing Sciences, Tokyo Institute of Technology, 2-12-1, Ookayama, Meguro-ku, Tokyo 152-8552, Japan\\
tera@is.titech.ac.jp
\vspace{.2cm}\\
Jun Ueki\\
Graduate School of Mathematical Science, The University of Tokyo,  3-8-1, Komaba, Meguro-ku, Tokyo, 153-8914, Japan\\
uekijun46@gmail.com}
\end{document}